\documentclass[a4paper,11pt,reqno]{article}
\usepackage{amsmath, amsfonts, amssymb, amsthm}
\input{epsf.sty}
\usepackage{mathrsfs}
\usepackage{amsmath}
\usepackage{amssymb}
\usepackage{graphicx}
\usepackage{subfigure}
\def \R {{\mathbb{R}}}
\numberwithin{equation}{section}

\begin{document}

\title{On the Cauchy problem for a generalized Camassa-Holm equation
with both quadratic and cubic nonlinearity}

\author{Xingxing $\mbox{Liu}^1$\footnote{E-mail: liuxingxing123456@163.com}
\quad Zhijun $\mbox{Qiao}^2$\footnote{E-mail:
qiao@utpa.edu} \quad and \quad Zhaoyang $\mbox{Yin}^1$\footnote{E-mail: mcsyzy@mail.sysu.edu.cn}\\
 $^1\mbox{Department}$ of Mathematics,
Sun Yat-sen University,\\ Guangzhou, Guangdong 510275 , China \\
$^2\mbox{Department}$ of Mathematics, University of Texas-Pan
American,\\ Edinburg, Texas 78541, USA\\}
\date{}
\maketitle

\begin{abstract}
In this paper, we study the Cauchy problem for a generalized
integrable Camassa-Holm equation with both quadratic and cubic
nonlinearity. By overcoming the difficulties caused by the
complicated mixed nonlinear structure, we firstly establish the
local well-posedness result in Besov spaces, and then present a
precise blow-up scenario for
strong solutions. Furthermore, we show the existence of single peakon by the method of analysis.\\

\noindent 2000 Mathematics Subject Classification: 35G25, 35L05
\smallskip\par
\noindent \textit{Keywords}: Generalized Camassa-Holm equation;
Besov spaces; local well-posedness; blow-up scenario; peakon.

\end{abstract}

\section{Introduction}
\par

In this paper, we consider the following partial differential
equation with both quadratic and cubic nonlinear terms
\begin{equation}\label{main}
m_t=\frac{1}{2}k_1\big((u^2-u_x^2)m\big)_x+\frac{1}{2}k_2(um_x+2m
u_x), \quad m=u-u_{xx},
\end{equation}
where $k_1,k_2$ are arbitrary constants. Eq. (\ref{main}) was first
proposed by Fokas in \cite{Fokas,Fokas1}. Very recently, it has been
shown that Eq. (\ref{main}) has a Lax pair, and can be written as
bi-Hamiltonian structure \cite{Qiao3}
\begin{equation*}
m_t=J\frac{\delta H_1}{\delta m}=K\frac{\delta H_2}{\delta m},
\end{equation*}
where
\begin{equation*}
J=k_1\partial m\partial^{-1}m\partial+\frac{1}{2}k_2(\partial m
+m\partial), \quad K=\partial-\partial^3,
\end{equation*}
and two Hamiltonians are
\begin{equation*}
H_1=\frac{1}{2}\int_{\R}(u^2+u^2_x)dx,
\end{equation*}
and 
 \begin{equation*}
H_2=\frac{1}{8}\int_{\R}(k_1u^4+2k_1u^2u^2_x-\frac{1}{3}k_1u^4_x+2k_2u^3+2k_2uu^2_x)dx.
\end{equation*}
Thus, Eq. (\ref{main}) is completely integrable.

Obviously, for $k_1=0,k_2=-2$, Eq. (\ref{main}) is reduced to the
Camassa-Holm (CH) equation \cite{Camassa1,FF}
\begin{equation*}
m_t+um_x+2u_xm =0,\quad m=u-u_{xx},
\end{equation*}
which describes the unidirectional propagation of waves at the free
surface of shallow water under the influence of gravity
\cite{Camassa1,Camassa2}. $u(t,x)$ stands for the fluid velocity at
time $t$ in the spatial $x$-direction, $x\in \R$, and $m(t,x)$
represents its potential density.

In the past few years, a large
amount literature has devoted to the investigation of the CH
equation, because it can describe both wave breaking phenomenon
\cite{Constantin3,Constantin5,Constantin6,Danchin,Li} (the solution
remains bounded while the slope of $u(t,x)$ becomes unbounded in
finite time), and solitary waves interacting like solitons
\cite{Camassa1,C-S,C-S3,C-S4,L-Y}. The well-posedness of the CH
equation has been shown in \cite{H-M,Li,Rodriguez} with the initial
data $u_0\in H^s(\R),s>\frac{3}{2}$ . In particular, Danchin
\cite{Danchin} has dealt with the initial-value problem of the CH
equation for the initial data in the Besov space $B_{p,r}^s,$ with
$1\leq p,r\leq +\infty,s>\max\{1+\frac{1}{p},\frac{3}{2}\}.$
However, the Cauchy problem of the CH equation is not locally
well-posed in $H^s(\R),s<\frac{3}{2}$. Indeed, the solution can not
depend uniformly continuously with respect to the initial data
\cite{H-M}. On the other hand, the CH equation has the peaked
solitions (peakons) of the form $\varphi_c(t,x)=ce^{-|x-ct|}$ with
the traveling speed $c>0$. For the peakon solution, we know that it
replicates a feature that is characteristic for the waves of great
height-waves of largest amplitude that are exact solutions of the
governing equations for water waves \cite{Con,Con1,Con-J,Toland}.
Constantin and Strauss \cite{C-S} gave an impressive proof of
stability of peakons by using the conservation laws.

For $k_1=-2,k_2=0$, Eq. (\ref{main}) becomes the following
equation with cubic nonlinearity
\begin{equation}\label{cubic}
m_t+bu_x+\big((u^2-u_x^2)m\big)_x=0,\quad m=u-u_{xx}, \quad b=const.
\end{equation}
which was derived independently by Fokas \cite{Fokas},
Fuchssteiner \cite{Fuchssteiner}, Olver and Rosenau \cite{Olver},
and Qiao \cite{Qiao1}. Eq. (\ref{cubic}) regains attention
due to its cuspon and peakon solution property and Lax pair \cite{Qiao1}, which may allow the initial value problem of
(\ref{cubic}) to be solved by the inverse scattering transform (IST)
method.
Unlike the CH equation, Eq. (\ref{cubic})
admits not only new cusp solitons (cuspons), but also possesses weak
kink solutions
 ($u,u_x,u_t$ are continuous, but $u_{xx}$
has a jump at its peak point) \cite{Qiao4,Qiao3}. It also has
significant differences from the CH equation about the dynamics of
the two-peakons and peakon-kink solutions \cite{Qiao3}. Recently,
the so called ''white'' solitons and ''dark'' ones of Eq.
(\ref{cubic}) have been presented in \cite{Sakovich} and \cite{I-L},
respectively. In \cite{Bies}, the authors apply the geometric and
analytic approaches to give a geometric interpretation to the
variable $m(t,x)$ and construct an infinite-dimensional Lie algebra
of symmetries to Eq. (\ref{cubic}). In \cite{G-L-O-Q}, the authors
consider the formulation of the singularities of solutions and show
that some solutions with certain initial date will blow up in finite
time, then they discuss the existence of single peakon of the form
$\varphi_c(t,x)=\pm\sqrt{\frac{3c}{2}}e^{-|x-ct|},c>0$, and
multi-peakon solutions for Eq. (\ref{cubic}). Very recently, the
orbital stability of peakons for Eq. (\ref{cubic}) has been proven
in \cite{Q-L-L}.

In the present paper, motivated by the study of the CH equation
\cite{Danchin}, our main work is to prove the local well-posedness
to the Cauchy problem (\ref{main}) in the nonhomogeneous Besov
spaces. However, one of the differences with \cite{Danchin} is that
we are required to deal with cubic nonlinearity in Besov spaces.
Moreover, the nonlinear term $"m_xu^2_x"$ makes us have to solve a
transport equation satisfied by $m$, rather than $u$. In contrast to
the case of the CH equation with initial data $u_0$ in the Sobolev
space $H^s(\R),s>\frac{3}{2}$, we can only prove the well-posedness
result with the initial profile $u_0$ in $H^s(\R),s>\frac{5}{2}$. In
our procedure, we have overcome the critical index case by the
interpolation method when we applied the transport theory to Eq.
(\ref{main}). Another one of the differences with \cite{Danchin,fu}
is that Eq. (\ref{main}) possesses the complicated mixed structure
nonlinear structure (with both quadratic and cubic nonlinearity). To
get the uniform boundedness of the approximate solutions
$\{u^{(n)}\}_{n\in \mathbb{N}}$, we have to handle the quadratic and
cubic nonlinear terms together in Eq. (\ref{main}). To overcome
these difficulties, we need to consider two cases: the small initial
data and the large one. Then with the local well-posedness result,
we may naturally present a precise blow-up scenario to Eq.
(\ref{main}) by combining the blow-up criterion of the CH equation
and the one of Eq. (\ref{cubic}).

The entire paper is organized as follows. In Section 2, we present
some facts on Besov spaces, some preliminary properties and the
transport equation theory. In Section 3, we establish the local
well-posedness result of Eq. (\ref{main}) in Besov spaces. In
Section 4, we derive a blow-up scenario for strong solutions to Eq.
(\ref{main}). In Section 5, we show that the existence of peakons
which can be understand as weak solutions for Eq. (\ref{main}).

$Notation.$ In the following, we denote $C>0$ a generic constant
only depending on $p,r,s$. Since our discussion about Eq.
(\ref{main}) is mainly on the line $\R$, for simplicity, we omit
$\R$ in our notations of function spaces. And we denote the Fourier
transform of a function $u$ as $\mathcal{F}u$.

\section{Preliminaries}
\newtheorem {remark2}{Remark}[section]
\newtheorem {proposition2}{Proposition}[section]
\newtheorem {definition2}{Definition}[section]
\newtheorem{theorem2}{Theorem}[section]
\newtheorem{lemma2}{Lemma}[section]
In this section, we will recall some basic theory of the
Littlewood-Paley decomposition and the transport equation theory on
Besov spaces, which will play an important role in the sequel. One
may get more details from \cite{Bahouri,Danchin}.

\begin{proposition2}\cite{Bahouri}
(Littlewood-Paley decomposition) Let
$\mathcal{B}:=\{\xi\in\mathbb{R},|\xi|\leq\frac{4}{3}\}$ and
$\mathcal{C}:=\{\xi\in\mathbb{R},\frac{3}{4}\leq|\xi|\leq\frac{8}{3}\}$.
Then there exist $\psi(\xi)\in C^{\infty}_c(\mathcal{B})$ and
$\varphi(\xi)\in C^{\infty}_c(\mathcal{C})$ such that
\begin{equation*}
\psi(\xi)+\sum_{q\in\mathbb{N}}\varphi(2^{-q}\xi)=1,\quad \forall\
\xi\in\mathbb{R}.
\end{equation*}
and
\begin{equation*}
Supp \varphi(2^{-q}\cdot) \cap Supp \varphi(2^{-q'}\cdot)
=\emptyset, \quad \mbox{if} \ |q-q'|\geq 2,
\end{equation*}
\begin{equation*}
Supp\psi(\cdot) \cap Supp \varphi(2^{-q}\cdot) =\emptyset, \quad
\mbox{if} \ q\geq 1.
\end{equation*}
Then for all $u\in \mathcal{S'}$($\mathcal{S'}$ denotes the tempered
distribution spaces), we can define the nonhomogeneous
Littlewood-Paley decomposition of a distribution $u$.
$$u=\sum_{q\in \mathbb{Z}}\Delta_qu,$$
where the localization operators are defined as follows:
\begin{equation*}
\Delta_qu:=0, \quad for \ q\leq-2,\quad
\Delta_{-1}u:=\psi(D)u=\mathcal{F}^{-1}(\psi\mathcal{F}u),
\end{equation*}
and
\begin{equation*}
\Delta_qu:=\varphi(2^{-q}D)u=\mathcal{F}^{-1}(\varphi(2^{-q}\xi)\mathcal{F}u),\quad
\mbox{for}\ q\geq0.
\end{equation*}
Furthermore, we can define the low frequency cut-off operator $S_q$
as follows:
$$S_q u:=\sum^{q-1}_{i=-1}\Delta_i u=\psi(2^{-q}D)u=\mathcal{F}^{-1}(\psi(2^{-q}\xi)\mathcal{F}u).$$
\end{proposition2}

\begin{definition2}\cite{Bahouri}
 (Besov spaces) Let $s\in \R,\ 1\leq p,r\leq\infty$. The nonhomogeneous Besov
 space $B^s_{p,r}(\mathbb{R})$ ($B^s_{p,r}$ for short) is defined by
 $$B^s_{p,r}:=\{u\in\mathcal{S'}(\mathbb{R}); \ \|u\|_{B^s_{p,r}}<\infty\},$$
 where
 \[\|u\|_{B^s_{p,r}}:=\left\{
\begin{array}{ll}
(\sum\limits_{q\in \mathbb{Z}}2^{qsr}\|\Delta_qu\|_{L^p}^r)^{\frac{1}{r}},  &{r<\infty,}\\
\sup\limits_{q\in \mathbb{Z}} 2^{qs}\|\Delta_qu\|_{L^p},
&{r=\infty.}
\end{array}
\right.
\]
If $s=\infty$, $B^\infty_{p,r}:=\bigcap\limits_{s\in \R}B^s_{p,r}.$
\end{definition2}

In order to state the local well-posedness result, we need to define
the following spaces.
\begin{definition2}
Let $T>0,s\in \R$ and $1\leq p\leq \infty$. We define
\begin{equation*}
E^s_{p,r}(T):=C([0,T];B^s_{p,r})\cap C^1([0,T];B^{s-1}_{p,r}), \quad
\mbox{for} \ r<\infty,
\end{equation*}
\begin{equation*}
E^s_{p,\infty}(T):=L^\infty([0,T];B^s_{p,\infty})\cap
Lip([0,T];B^{s-1}_{p,\infty}),
\end{equation*}
and
$$E^s_{p,r}:=\bigcap_{T>0}E^s_{p,r}(T).$$
\end{definition2}

 Next, we list the following useful properties for Besov spaces.
 \begin{proposition2}\cite{Bahouri,Danchin}
Let $s\in \R,\ 1\leq p,r,p_i,r_i\leq\infty,i=1,2$. Then\\
(i) Density: if $1\leq p,r<\infty$, then $\mathcal{C}^\infty_c$ is
dense in $B^s_{p,r}$.\\
(ii) Embedding: $B^s_{p_1,r_1}\hookrightarrow
B^{s-(\frac{1}{p_1}-\frac{1}{p_2})}_{p_2,r_2}, \quad \mbox{for} \
p_1\leq p_2\ \mbox{and} \ r_1\leq r_2$,
\begin{equation*}
B^{s_2}_{p,r_2}\hookrightarrow B^{s_1}_{p,r_1} \ \mbox{locally
compact}, \quad \mbox{for} \ s_1<s_2.
\end{equation*}\\
(iii) Algebraic properties: if $s>0$, $B^s_{p,r}\cap L^\infty$ is an
algebra. Furthermore, $B^s_{p,r}$ is an algebra, provided that
$s>\frac{1}{p}$ or $s\geq\frac{1}{p} \ \mbox{and}\  r=1.$\\
(iv) Fatou lemma: if $\{u^{(n)}\}_{n\in \mathbb{N}}$ is bounded in
$B^s_{p,r}$ and tends to $u$ in $\mathcal{S'}$, then $u\in
B^s_{p,r}$. Moreover,
$$\|u\|_{B^s_{p,r}}\leq \lim \inf_{n\rightarrow
\infty}\|u^{(n)}\|_{B^s_{p,r}}.$$
(v) Complex interpolation: if
$u\in B^{s_1}_{p,r}\cap B^{s_2}_{p,r}$, then for all
$\theta\in[0,1]$, we have $u\in B^{\theta s_1+(1-\theta)s_2}_{p,r}$.
Moreover,
$$\|u\|_{B^{\theta s_1+(1-\theta)s_2}_{p,r}}\leq
\|u\|^\theta_{B^{s_1}_{p,r}}\|u\|^{(1-\theta)}_{B^{s_2}_{p,r}}.$$
(vi) One-dimensional Morse-type estimate:\\
1) If $s>0$,
\begin{equation*}
\|uv\|_{B^s_{p,r}}\leq C
(\|u\|_{B^s_{p,r}}\|v\|_{L^\infty}+\|u\|_{L^\infty}\|v\|_{B^s_{p,r}}).
\end{equation*}
2) If $s_1\leq \frac{1}{p}$
 , $s_2>\frac{1}{p}$($s_2\geq\frac{1}{p}\ \mbox{if}\ r=1$) and $s_1+s_2>0,$
 \begin{eqnarray*}
\|u v\|_{B^{s_1}_{p,r}}\leq C
\|u\|_{B^{s_1}_{p,r}}\|v\|_{B^{s_2}_{p,r}},
\end{eqnarray*}
where $C$ is a constant independent of $u$ and $v$.\\
(vii) Action of Fourier multipliers on Besov spaces: let $m\in \R$
and $f$ be a $S^m$-multiplier ($i.e.,$ $f:\R\rightarrow \R$ is a
smooth function and satisfies that for each multi-index $\alpha$,
there exists a constant $C_\alpha$ such that $|\partial^\alpha
f(\xi)|\leq C_\alpha(1+|\xi|)^{m-|\alpha|}$, for $\forall \xi \in
\R.$) Then the operator $f(D)$ is continuous from $B^s_{p,r}$ to
$B^{s-m}_{p,r}$.
\end{proposition2}

Now we state the following transport equation theory that is crucial
to prove local well-posedness for Eq. (\ref{main}).
\begin{lemma2}\cite{Bahouri,Danchin} (A priori estimate)
Let $1\leq p,r\leq +\infty$ and
$s>-\min\{\frac{1}{p},1-\frac{1}{p}\}$. Assume that $v$ be a
function such that
 $\partial_x v$ belongs to $L^1([0,T];B^{s-1}_{p,r})$ if $s>1+\frac{1}{p}$
 or to $L^1([0,T];B^{\frac{1}{p}}_{p,r}\cap L^\infty)$ otherwise.
Suppose also that $f_0\in B^s_{p,r}$, $F\in L^1([0,T];B^s_{p,r})$,
and that $f\in L^\infty([0,T];B^{s}_{p,r})\cap
C([0,T];\mathcal{S'})$
 be the solution of the one-dimensional transport equation
 \begin{equation}
\left\{\begin{array}{ll}\partial_tf+v\cdot\partial_x f=F,\\
 f\mid_{t=0}=f_0.\\
\end{array}\right.
\end{equation}
Then there exists a constant $C$ depending only on $s,p,r$ such that
the following statements hold for $t\in [0,T]$\\
(i) If $r=1$ or
$s\neq 1+\frac{1}{p}$,
\begin{eqnarray*}
\|f(t)\|_{B^{s}_{p,r}}\leq
\|f_0\|_{B^{s}_{p,r}}+\int_0^t\|F(\tau)\|_{B^{s}_{p,r}}d\tau
+C\int_0^tV'(\tau)\|f(\tau)\|_{B^{s}_{p,r}}d\tau,
\end{eqnarray*}
 or
 \begin{eqnarray*}
\|f(t)\|_{B^{s}_{p,r}}\leq
e^{CV(t)}(\|f_0\|_{B^{s}_{p,r}}+\int_0^te^{-CV(\tau)}
\|F(\tau)\|_{B^{s}_{p,r}}d\tau),
\end{eqnarray*}
where \[V(t)=\left\{
\begin{array}{ll}
\int_0^t\|\partial_xv(\tau,\cdot)\|_{(B^{\frac{1}{p}}_{p,r}\cap L^\infty)}d\tau, &{s<1+\frac{1}{p},}\\
\int_0^t\|\partial_xv(\tau,\cdot)\|_{B^{s-1}_{p,r}}d\tau,
&{s>1+\frac{1}{p}.}
\end{array}
\right.
\]
(ii) If $s\leq 1+\frac{1}{p}$, and $\partial_xf_0,\partial_xf\in
L^\infty([0,T]\times \R)$ and $\partial_xF\in L^1([0,T];L^\infty)$,
then
\begin{eqnarray*}
&&\|f(t)\|_{B^{s}_{p,r}}+\|\partial_x f(t)\|_{L^\infty}\\
&\leq& e^{CV(t)}(\|f_0\|_{B^{s}_{p,r}}+\|\partial_x
f_0\|_{L^\infty}+ \int_0^te^{-CV(\tau)}(
\|F(\tau)\|_{B^{s}_{p,r}}+\|\partial_xF(\tau)\|_{L^\infty})d\tau),
\end{eqnarray*}
where
$V(t)=\int_0^t\|\partial_xv(\tau,\cdot)\|_{(B^{\frac{1}{p}}_{p,r}\cap
L^\infty)}d\tau.$\\
(iii) If $f=v$, then for all $s>0$, the estimate in (i) holds with
$V(t)=\int_0^t\|\partial_xv(\tau,\cdot)\|_{
L^\infty}d\tau.$\\
(iv) If $r<\infty$, then $f\in C([0,T];B^{s}_{p,r})$. If $r=\infty$,
then $f\in C([0,T];B^{s'}_{p,1})$ for all $s'<s$.
\end{lemma2}

\begin{lemma2}\cite{Danchin}
(Existence and uniqueness) Let $p,r,s,f_0$ and $F$ be as in the
statement of Lemma 2.1. Suppose that $v\in
L^\rho([0,T];B_{\infty,\infty}^{-M})$ for some $\rho>1,M>0$ and
$\partial_x v\in L^1([0,T];B^{\frac{1}{p}}_{p,\infty}\cap L^\infty)$
if $s<1+\frac{1}{p}$, and $\partial_x v\in L^1([0,T];B^{s-1}_{p,r})$
if $s>1+\frac{1}{p}$ or $s=1+\frac{1}{p}$ and $r=1$. Then the
transport equation (2.1) has a unique solution $f\in
L^\infty([0,T];B^{s}_{p,r})\bigcap
(\bigcap\limits_{s'<s}C([0,T];B^{s'}_{p,1}))$ and the corresponding
inequalities in Lemma 2.1 hold true. Moreover, if $r<\infty$, then
$f\in C([0,T];B^{s}_{p,r}).$
\end{lemma2}

\section{Local well-posedness}
\newtheorem {remark3}{Remark}[section]
\newtheorem{theorem3}{Theorem}[section]
\newtheorem{lemma3}{Lemma}[section]
\par
In this section, we shall study the local well-posedness of Eq.
(\ref{main}) in the nonhomogeneous Besov spaces. At first, we
present a priori estimates about the solutions of Eq. (\ref{main}),
which can be applied to prove the uniqueness and continuity with the
initial data in some sense.

\begin{lemma3}
Suppose that $1\leq p,r\leq\infty$ and
$s>\max\{2+\frac{1}{p},\frac{5}{2},3-\frac{1}{p}\}$. Let
$u^{(1)},u^{(2)}\in L^\infty([0,T];B^{s}_{p,r})\cap
C([0,T];\mathcal{S'})$ be two given solutions to Eq. (\ref{main})
with initial data $u_0^{(1)},u_0^{(2)}\in B^{s}_{p,r}$, and let
$u^{(12)}:=u^{(2)}-u^{(1)}$ and $m^{(12)}:=m^{(2)}-m^{(1)}$. Then
for all $t\in [0,T]$, we have\\
(1) If $s>\max\{2+\frac{1}{p},\frac{5}{2},3-\frac{1}{p}\}$ and
$s\neq 4+\frac{1}{p}$, then
\begin{eqnarray*}
\|u^{(12)}\|_{B^{s-1}_{p,r}}&\leq&
\|u_0^{(12)}\|_{B^{s-1}_{p,r}}\exp
\{C\int^t_0(\|u^{(1)}(\tau)\|^2_{B^{s}_{p,r}}
+\|u^{(2)}(\tau)\|^2_{B^{s}_{p,r}}\\
&
&+\|u^{(1)}(\tau)\|_{B^{s}_{p,r}}
+\|u^{(2)}(\tau)\|_{B^{s}_{p,r}})d\tau\}.
\end{eqnarray*}
(2) If $s= 4+\frac{1}{p}$, then
\begin{eqnarray*}
&&\|u^{(12)}\|_{B^{s-1}_{p,r}}\\
&\leq& C
\|u_0^{(12)}\|^{\theta}_{B^{s-1}_{p,r}}(\|u^{(1)}\|_{B^{s}_{p,r}}+\|u^{(2)}\|_{B^{s}_{p,r}})^{1-\theta}
\exp\{\theta C\int^t_0(\|u^{(1)}(\tau)\|^2_{B^{s}_{p,r}}\\
& &
+\|u^{(2)}(\tau)\|^2_{B^{s}_{p,r}}+\|u^{(1)}(\tau)\|_{B^{s}_{p,r}}
+\|u^{(2)}(\tau)\|_{B^{s}_{p,r}})d\tau\},
\end{eqnarray*}
with $\theta=\frac{1}{2}(1-\frac{1}{2p})\in(0,1)$.
\end{lemma3}
\begin{proof}
It is obvious that $u^{(12)}\in L^\infty([0,T];B^{s}_{p,r})\cap
C([0,T];\mathcal{S'})$ and $u^{(12)},m^{(12)}$ solves the following
transport equation
\begin{equation}\ \
\left\{\begin{array}{ll}m^{(12)}_t-\{\frac{k_1}{2}[(u^{(1)})^2-(u_x^{(1)})^2]
+\frac{k_2}{2}u^{(1)}\}m_x^{(12)}
=F(t,x), & t>0,x\in \mathbb{R},\\
m^{(12)}\mid_{t=0}=m_0^{(12)}:=m_0^{(2)}-m_0^{(1)}, \quad
m=u-u_{xx},& x\in \mathbb{R},
\end{array}\right.
\end{equation}
where
$F(t,x):=\frac{k_1}{2}[u^{(12)}(u^{(1)}+u^{(2)})-u_x^{(12)}(u_x^{(1)}+u_x^{(2)})]\partial_xm^{(2)}
+k_1[u_x^{(12)}(m^{(2)})^2+u_x^{(1)}m^{(12)}(m^{(1)}+m^{(2)})]+\frac{k_2}{2}u^{(12)}m_x^{(2)}
+k_2u_x^{(12)}m^{(2)}+k_2u_x^{(1)}m^{(12)}.$

Applying Lemma 2.1 to the transport equation (3.1), we have
\begin{eqnarray}
\|m^{(12)}\|_{B^{s-3}_{p,r}}&\leq&
\|m_0^{(12)}\|_{B^{s-3}_{p,r}}+\int^t_0
\|F(\tau)\|_{B^{s-3}_{p,r}}d\tau \nonumber\\
&
&+C\int^t_0\|(u^{(1)})^2-(u_x^{(1)})^2+u^{(1)}\|_{B^{s-3}_{p,r}}\|m^{(12)}\|_{B^{s-3}_{p,r}}d\tau.
\end{eqnarray}
Indeed, if $\max\{2+\frac{1}{p},\frac{5}{2}\}<s\leq 3+\frac{1}{p}$,
by Proposition 2.2 (vi), we get
\begin{eqnarray*}
&&\|F(\tau)\|_{B^{s-3}_{p,r}}\\
 &\leq& C \{
\|[u^{(12)}(u^{(1)}+u^{(2)})-u_x^{(12)}(u_x^{(1)}+u_x^{(2)})]\partial_xm^{(2)}\|_{B^{s-3}_{p,r}}\\
&&+\|u_x^{(12)}(m^{(2)})^2\|_{B^{s-3}_{p,r}}
+\|u_x^{(1)}m^{(12)}(m^{(1)} +m^{(2)})\|_{B^{s-3}_{p,r}}\\&&
+\|u^{(12)}m_x^{(2)}+u_x^{(12)}m^{(2)}+u_x^{(1)}m^{(12)}\|_{B^{s-3}_{p,r}}
\}\\
&\leq& C
\{\|m_x^{(2)}\|_{B^{s-3}_{p,r}}\|[u^{(12)}(u^{(1)}+u^{(2)})-u_x^{(12)}(u_x^{(1)}+u_x^{(2)})]\|_{B^{s-2}_{p,r}}\\
&&+\|u_x^{(12)}\|_{B^{s-3}_{p,r}}\|(m^{(2)})^2\|_{B^{s-2}_{p,r}}+\|m^{(12)}\|_{B^{s-3}_{p,r}}
\|u_x^{(1)}
(m^{(1)} +m^{(2)})\|_{B^{s-2}_{p,r}}\\
&&+\|m_x^{(2)}\|_{B^{s-3}_{p,r}}\|u^{(12)}\|_{B^{s-2}_{p,r}}+\|u_x^{(12)}\|_{B^{s-3}_{p,r}}\|m^{(2)}\|_{B^{s-2}_{p,r}}\\
&&+\|m^{(12)}\|_{B^{s-3}_{p,r}}\|u_x^{(1)}\|_{B^{s-2}_{p,r}} \}.
\end{eqnarray*}
Since $s>\max\{2+\frac{1}{p},\frac{5}{2}\}$, we know that
${B^{s-2}_{p,r}}$ is an algebra. Thus, we deduce
\begin{eqnarray}
&&\|F(\tau)\|_{B^{s-3}_{p,r}}\nonumber\\
&\leq&
C[\|u^{(12)}\|_{B^{s-1}_{p,r}}(\|u^{(1)}\|^2_{B^{s}_{p,r}}+\|u^{(2)}\|^2_{B^{s}_{p,r}}
+\|u^{(1)}\|_{B^{s}_{p,r}}+\|u^{(2)}\|_{B^{s}_{p,r}}).
\end{eqnarray}
For $s>3+\frac{1}{p}$, the inequality (3.3) also holds true in view
of the fact that ${B^{s-3}_{p,r}}$ is an algebra.\\
Note that
\begin{eqnarray*}
&&\|(u^{(1)})^2-(u_x^{(1)})^2+u^{(1)}\|_{B^{s-3}_{p,r}}\|m^{(12)}\|_{B^{s-3}_{p,r}}\\
&\leq&
C[\|u^{(12)}\|_{B^{s-1}_{p,r}}(\|u^{(1)}\|^2_{B^{s}_{p,r}}+\|u^{(2)}\|^2_{B^{s}_{p,r}}
+\|u^{(1)}\|_{B^{s}_{p,r}}+\|u^{(2)}\|_{B^{s}_{p,r}}).\end{eqnarray*}
Therefore, inserting the above inequality and (3.3) into (3.2), we
obtain
\begin{eqnarray*}
\|u^{(12)}\|_{B^{s-1}_{p,r}}&\leq&
\|u_0^{(12)}\|_{B^{s-1}_{p,r}}+C\int^t_0
[\|u^{(12)}(\tau)\|_{B^{s-1}_{p,r}}\\
&& \times (\|u^{(1)}\|^2_{B^{s}_{p,r}}+\|u^{(2)}\|^2_{B^{s}_{p,r}}
+\|u^{(1)}\|_{B^{s}_{p,r}}+\|u^{(2)}\|_{B^{s}_{p,r}})] d\tau.
\end{eqnarray*}
Then, by Gronwall's inequality, we prove (1).

Since we can not apply Lemma 2.1 to (3.1) for the critical case $s=
4+\frac{1}{p}$, we here use the interpolation method to deal with
it.

In fact, we can choose $\theta=\frac{1}{2}(1-\frac{1}{2p})\in(0,1)$,
such that $s-1=3+\frac{1}{p}=(1-\theta)(4+\frac{1}{2p})+\theta
(2+\frac{1}{2p})$. Then, by Proposition 2.2 (v), we have
$$\|u^{(12)}\|_{B^{s-1}_{p,r}}=\|u^{(12)}\|_{B^{3+\frac{1}{p}}_{p,r}}
\leq \|u^{(12)}\|^\theta_{B^{2+\frac{1}{2p}}_{p,r}}
\|u^{(12)}\|^{1-\theta}_{B^{4+\frac{1}{2p}}_{p,r}}.$$ Then, from the
obtained result of (1) in this lemma, we get
\begin{eqnarray*}
\|u^{(12)}\|_{B^{s-1}_{p,r}} &\leq&
[\|u_0^{(12)}\|_{B^{2+\frac{1}{2p}}_{p,r}} \exp(
C\int^t_0(\|u^{(1)}\|^2_{B^{3+\frac{1}{2p}}_{p,r}}+\|u^{(2)}\|^2_{B^{3+\frac{1}{2p}}_{p,r}}\\
&&+\|u^{(1)}\|_{B^{3+\frac{1}{2p}}_{p,r}}+\|u^{(2)}\|_{B^{3+\frac{1}{2p}}_{p,r}})
d\tau) ]^\theta
\times(\|u^{(1)}\|_{B^{4+\frac{1}{2p}}_{p,r}}+\|u^{(2)}\|_{B^{4+\frac{1}{2p}}_{p,r}})^{1-\theta}\\
&\leq& C
\|u_0^{(12)}\|^{\theta}_{B^{s-1}_{p,r}}(\|u^{(1)}\|_{B^{s}_{p,r}}+\|u^{(2)}\|_{B^{s}_{p,r}})^{1-\theta}
\exp\{\theta C\int^t_0(\|u^{(1)}(\tau)\|^2_{B^{s}_{p,r}}\\
& &
+\|u^{(2)}(\tau)\|^2_{B^{s}_{p,r}}+\|u^{(1)}(\tau)\|_{B^{s}_{p,r}}
+\|u^{(2)}(\tau)\|_{B^{s}_{p,r}})d\tau\}.
\end{eqnarray*}
This completes the proof of Lemma 3.1.
\end{proof}

Next, we use the classical Friedrichs regularization method to
construct the approximation solutions to Eq. (\ref{main}).
\begin{lemma3}
Suppose that $p,r$ and $s$ be as in the statement of Lemma 3.1,
$u_0\in B^{s}_{p,r}$ and $u^{(0)}:=0$. Then there exists a sequence
of smooth functions $\{u^{(n)}\}_{n\in \mathbb{N}} \in
C(\R^+;B^{\infty}_{p,r})$ solving the following transport equation
by induction
\begin{equation}
\left\{\begin{array}{ll}m^{(n+1)}_t-\{\frac{k_1}{2}[(u^{(n)})^2-(u_x^{(n)})^2]+\frac{k_2}{2}u^{(n)}\}m^{(n+1)}_x
=\\
\quad \quad \quad \quad \quad \quad \quad \quad \quad\quad \quad
\quad \quad
 k_1u_x^{(n)}(m^{(n)})^2+k_2u_x^{(n)}m^{(n)},  & t>0,x\in \mathbb{R}, \\
u^{(n+1)}\mid_{t=0}=u^{(n+1)}_{0}(x)=S_{n+1}u_0, & x\in \mathbb{R}.
\end{array}\right.
\end{equation}
Moreover, there exists a $T>0$ such that the solutions satisfying
the following properties.\\
(1) $\{u^{(n)}\}_{n\in \mathbb{N}}$ is uniformly bounded in $E^{s}_{p,r}(T)$.\\
(2) $\{u^{(n)}\}_{n\in \mathbb{N}}$ is a Cauchy sequence in
$C([0,T];B^{s-1}_{p,r}).$
\end{lemma3}

\begin{proof}By the definition of $S_q$, we know that all the data $S_{n+1}u_0 \in
B^{\infty}_{p,r}$. Thus, from Lemma 2.2, we deduce by induction that
for all $n\in \mathbb{N}$, Eq. (3.4) has a global solution belonging
$C(\R^+;B^{\infty}_{p,r})$.

For $s>\max\{2+\frac{1}{p},\frac{5}{2},3-\frac{1}{p}\}$ and $s\neq
4+\frac{1}{p}$, by Lemma 2.1, we obtain
\begin{eqnarray}
\|m^{(n+1)}\|_{B^{s-2}_{p,r}} &\leq&
e^{C\int^t_0\|[(u^{(n)})^2-(u_x^{(n)})^2+u^{(n)}](\tau)\|_{B^{s-2}_{p,r}}d\tau}(\|m_0\|_{B^{s-2}_{p,r}}\nonumber\\
&&+C\int^t_0e^{-C\int^\tau_0\|[(u^{(n)})^2-(u_x^{(n)})^2+u^{(n)}](\tau')\|_{B^{s-2}_{p,r}}d\tau'}\nonumber\\
&& \times
\|u_x^{(n)}(m^{(n)})^2+u_x^{(n)}m^{(n)}\|_{B^{s-2}_{p,r}}d\tau).
\end{eqnarray}
Since $s>2+\frac{1}{p}$, we know that $B^{s-2}_{p,r}$ is an algebra.
Thus, we have
\begin{eqnarray*}
\|u_x^{(n)}(m^{(n)})^2+u_x^{(n)}m^{(n)}\|_{B^{s-2}_{p,r}}
 &\leq&
C\|u_x^{(n)}\|_{B^{s-2}_{p,r}}(\|(m^{(n)})^2\|_{B^{s-2}_{p,r}}+\|(m^{(n)})\|_{B^{s-2}_{p,r}})\\
&\leq& C(\|u^{(n)}\|^3_{B^{s}_{p,r}}+\|u^{(n)}\|^2_{B^{s}_{p,r}}),
\end{eqnarray*}
and
\begin{eqnarray*}
\|(u^{(n)})^2-(u_x^{(n)})^2+u^{(n)}\|_{B^{s-2}_{p,r}} \leq
C(\|u^{(n)}\|^2_{B^{s}_{p,r}}+\|u^{(n)}\|_{B^{s}_{p,r}}).
\end{eqnarray*}
Inserting the above inequalities into (3.5), we obtain
\begin{eqnarray}
\|u^{(n+1)}\|_{B^{s}_{p,r}} &\leq&
e^{C\int^t_0(\|u^{(n)}\|^2_{B^{s}_{p,r}}+\|u^{(n)}\|_{B^{s}_{p,r}})(\tau)d\tau}\|u_0\|_{B^{s}_{p,r}}\nonumber\\
&&+C\int^t_0e^{C\int_\tau^t(\|u^{(n)}\|^2_{B^{s}_{p,r}}+\|u^{(n)}\|_{B^{s}_{p,r}})(\tau')d\tau'}\nonumber\\
 &&\times (\|u^{(n)}\|^3_{B^{s}_{p,r}}+\|u^{(n)}\|^2_{B^{s}_{p,r}})d\tau.
\end{eqnarray}

In order to prove the uniform boundedness of $\{u^{(n)}\}_{n\in
\mathbb{N}}$, we shall divide our discussion into two parts. When
$2\|u_0\|_{B^{s}_{p,r}}<1$, we choose a $T_1>0$ such that
$$T_1<\min\{\frac{1-2\|u_0\|_{B^{s}_{p,r}}}{8C\|u_0\|_{B^{s}_{p,r}}},\frac{1}{4C}\},$$
and suppose by induction that for all $t\in[0,T_1]$
\begin{eqnarray}
\|u^{(n)}\|_{B^{s}_{p,r}}\leq
\frac{2\|u_0\|_{B^{s}_{p,r}}}{1-8C\|u_0\|_{B^{s}_{p,r}}t}.
\end{eqnarray}
Noting that $t\leq T_1<
\frac{1-2\|u_0\|_{B^{s}_{p,r}}}{8C\|u_0\|_{B^{s}_{p,r}}}$, we have
$\frac{2\|u_0\|_{B^{s}_{p,r}}}{1-8C\|u_0\|_{B^{s}_{p,r}}t}<1.$ By
(3.7), we obtain
\begin{eqnarray}
&&\|u^{(n)}\|_{B^{s}_{p,r}}\leq
\frac{2\|u_0\|_{B^{s}_{p,r}}}{1-8C\|u_0\|_{B^{s}_{p,r}}t}\nonumber\\
&\leq&
(\frac{2\|u_0\|_{B^{s}_{p,r}}}{1-8C\|u_0\|_{B^{s}_{p,r}}t})^{\frac{2}{3}}
 \leq
 (\frac{2\|u_0\|_{B^{s}_{p,r}}}{1-8C\|u_0\|_{B^{s}_{p,r}}t})^{\frac{1}{2}}.
\end{eqnarray}
From (3.8), we can deduce that
\begin{eqnarray*}
&&C\int_\tau^t(\|u^{(n)}\|^2_{B^{s}_{p,r}}+\|u^{(n)}\|_{B^{s}_{p,r}})(\tau')d\tau'\\
&\leq& C
\int_\tau^t\{[(\frac{2\|u_0\|_{B^{s}_{p,r}}}{1-8C\|u_0\|_{B^{s}_{p,r}}\tau'})^{\frac{1}{2}}]^2+
\frac{2\|u_0\|_{B^{s}_{p,r}}}{1-8C\|u_0\|_{B^{s}_{p,r}}\tau'}\}d\tau'\\
 &\leq&
-\frac{1}{2}\int_\tau^t
\frac{-8C\|u_0\|_{B^{s}_{p,r}}}{1-8C\|u_0\|_{B^{s}_{p,r}}\tau'}d\tau'\\
&=& \ln \sqrt{1-8C\|u_0\|_{B^{s}_{p,r}}\tau}-\ln
\sqrt{1-8C\|u_0\|_{B^{s}_{p,r}}t}.
\end{eqnarray*}
Inserting the above inequality and (3.8) into (3.6) yields
\begin{eqnarray}
&&\|u^{(n+1)}\|_{B^{s}_{p,r}}\\
&\leq&
\frac{\|u_0\|_{B^{s}_{p,r}}}{\sqrt{1-8C\|u_0\|_{B^{s}_{p,r}}t}} +
\frac{C}{\sqrt{1-8C\|u_0\|_{B^{s}_{p,r}}t}}\{\int^t_0\sqrt{1-8C\|u_0\|_{B^{s}_{p,r}}\tau}\nonumber\\
&&\times[((\frac{2\|u_0\|_{B^{s}_{p,r}}}{1-8C\|u_0\|_{B^{s}_{p,r}}\tau})^{\frac{2}{3}})^{3}
+((\frac{2\|u_0\|_{B^{s}_{p,r}}}{1-8C\|u_0\|_{B^{s}_{p,r}}\tau})^{\frac{1}{2}})^2]\}d\tau\nonumber\\
&\leq&
\frac{\|u_0\|_{B^{s}_{p,r}}}{\sqrt{1-8C\|u_0\|_{B^{s}_{p,r}}t}}
+\frac{\|u_0\|_{B^{s}_{p,r}}}{\sqrt{1-8C\|u_0\|_{B^{s}_{p,r}}t}}\int^t_0
\frac{4C\|u_0\|_{B^{s}_{p,r}}}{(1-8C\|u_0\|_{B^{s}_{p,r}}\tau)^{\frac{3}{2}}}d\tau\nonumber\\
&&+\frac{1}{\sqrt{1-8C\|u_0\|_{B^{s}_{p,r}}t}}\int^t_0
\frac{2C\|u_0\|_{B^{s}_{p,r}}}{(1-8C\|u_0\|_{B^{s}_{p,r}}\tau)^{\frac{1}{2}}}d\tau\nonumber\\
&\leq&\frac{\|u_0\|_{B^{s}_{p,r}}}{\sqrt{1-8C\|u_0\|_{B^{s}_{p,r}}t}}
+\frac{\|u_0\|_{B^{s}_{p,r}}}{\sqrt{1-8C\|u_0\|_{B^{s}_{p,r}}t}}
[(1-8C\|u_0\|_{B^{s}_{p,r}}\tau)^{-\frac{1}{2}}|^{t}_{0}]\nonumber\\
&&+\frac{1}{\sqrt{1-8C\|u_0\|_{B^{s}_{p,r}}t}}
[-\frac{1}{2}(1-8C\|u_0\|_{B^{s}_{p,r}}\tau)^{\frac{1}{2}}|^{t}_{0}]\nonumber\\
&\leq&\frac{\|u_0\|_{B^{s}_{p,r}}}{\sqrt{1-8C\|u_0\|_{B^{s}_{p,r}}t}}
(1+\frac{1}{\sqrt{1-8C\|u_0\|_{B^{s}_{p,r}}t}}-1)
+\frac{1}{\sqrt{1-8C\|u_0\|_{B^{s}_{p,r}}t}}\nonumber\\
&&\times\frac{1}{2}(1-\sqrt{1-8C\|u_0\|_{B^{s}_{p,r}}t})\nonumber\\
&\leq&
\frac{2\|u_0\|_{B^{s}_{p,r}}}{1-8C\|u_0\|_{B^{s}_{p,r}}t},\nonumber
\end{eqnarray}
where we used the following fact that $$T_1<\frac{1}{4C}\Rightarrow
\frac{1}{2}(1-\sqrt{1-8C\|u_0\|_{B^{s}_{p,r}}t})\leq
\frac{\|u_0\|_{B^{s}_{p,r}}}{\sqrt{1-8C\|u_0\|_{B^{s}_{p,r}}t}},$$
in the last inequality. Thus, we prove (3.7) for the case
$2\|u_0\|_{B^{s}_{p,r}}<1$.

On the other hand, when $2\|u_0\|_{B^{s}_{p,r}}\geq 1$, we choose a
$T_2>0$ such that $T_2\leq
\frac{1-e^{-1}}{16C\|u_0\|^2_{B^{s}_{p,r}}}<\frac{1}{16C\|u_0\|^2_{B^{s}_{p,r}}}$,
and suppose by induction that for all $t\in [0,T_2]$
\begin{eqnarray}
\|u^{(n)}\|_{B^{s}_{p,r}}\leq
\frac{2\|u_0\|_{B^{s}_{p,r}}}{\sqrt{1-16C\|u_0\|^2_{B^{s}_{p,r}}t}}.
\end{eqnarray}
Noting that $2\|u_0\|_{B^{s}_{p,r}}\geq 1$, we get
$\frac{2\|u_0\|_{B^{s}_{p,r}}}{(1-16C\|u_0\|^2_{B^{s}_{p,r}}t)^\frac{1}{2}}\geq1.$
From (3.10), we obtain
\begin{eqnarray}
&&\|u^{(n)}\|_{B^{s}_{p,r}}\leq
\frac{2\|u_0\|_{B^{s}_{p,r}}}{\sqrt{1-16C\|u_0\|^2_{B^{s}_{p,r}}t}}\nonumber\\
&\leq&
(\frac{2\|u_0\|_{B^{s}_{p,r}}}{\sqrt{1-16C\|u_0\|^2_{B^{s}_{p,r}}t}})^{\frac{3}{2}}
 \leq
 (\frac{2\|u_0\|_{B^{s}_{p,r}}}{\sqrt{1-16C\|u_0\|^2_{B^{s}_{p,r}}t}})^{2}.
\end{eqnarray}
By (3.11), we find that
\begin{eqnarray*}
&&C\int_\tau^t(\|u^{(n)}\|^2_{B^{s}_{p,r}}+\|u^{(n)}\|_{B^{s}_{p,r}})(\tau')d\tau'\\
&\leq& C
\int_\tau^t\{(\frac{2\|u_0\|_{B^{s}_{p,r}}}{\sqrt{1-16C\|u_0\|^2_{B^{s}_{p,r}}\tau'}})^2+
\frac{4\|u_0\|^2_{B^{s}_{p,r}}}{1-16C\|u_0\|^2_{B^{s}_{p,r}}\tau'}\}d\tau'\\
 &\leq&
-\frac{1}{2}\int_\tau^t
\frac{-16C\|u_0\|^2_{B^{s}_{p,r}}}{1-16C\|u_0\|^2_{B^{s}_{p,r}}\tau'}d\tau'\\
&=& \ln \sqrt{1-16C\|u_0\|^2_{B^{s}_{p,r}}\tau}-\ln
\sqrt{1-16C\|u_0\|^2_{B^{s}_{p,r}}t}.
\end{eqnarray*}
Inserting the above inequality and (3.11) into (3.6), we obtain
\begin{eqnarray}
&&\|u^{(n+1)}\|_{B^{s}_{p,r}}\\
&\leq&
\frac{\|u_0\|_{B^{s}_{p,r}}}{\sqrt{1-16C\|u_0\|^2_{B^{s}_{p,r}}t}} +
\frac{C}{\sqrt{1-16C\|u_0\|^2_{B^{s}_{p,r}}t}}\{\int^t_0\sqrt{1-16C\|u_0\|^2_{B^{s}_{p,r}}\tau}\nonumber\\
&&\times[(\frac{2\|u_0\|_{B^{s}_{p,r}}}{\sqrt{1-16C\|u_0\|^2_{B^{s}_{p,r}}\tau}})^{3}
+((\frac{2\|u_0\|_{B^{s}_{p,r}}}{\sqrt{1-16C\|u_0\|^2_{B^{s}_{p,r}}\tau}})^{\frac{3}{2}})^2]\}d\tau\nonumber\\
&\leq&
\frac{\|u_0\|_{B^{s}_{p,r}}}{\sqrt{1-16C\|u_0\|^2_{B^{s}_{p,r}}t}}
+\frac{\|u_0\|_{B^{s}_{p,r}}}{\sqrt{1-16C\|u_0\|^2_{B^{s}_{p,r}}t}}\int^t_0
\frac{16C\|u_0\|^2_{B^{s}_{p,r}}}{(1-16C\|u_0\|^2_{B^{s}_{p,r}}\tau)}d\tau\nonumber\\
&\leq&\frac{\|u_0\|_{B^{s}_{p,r}}}{\sqrt{1-16C\|u_0\|^2_{B^{s}_{p,r}}t}}
[1-\ln (1-16C\|u_0\|^2_{B^{s}_{p,r}}\tau)|^{t}_{0}]\nonumber\\
&\leq&\frac{2\|u_0\|_{B^{s}_{p,r}}}{\sqrt{1-16C\|u_0\|^2_{B^{s}_{p,r}}t}},\nonumber
\end{eqnarray}
where we used the following fact that $$T_2\leq
\frac{1-e^{-1}}{16C\|u_0\|^2_{B^{s}_{p,r}}}\Rightarrow 1-\ln
(1-16C\|u_0\|^2_{B^{s}_{p,r}}t)\leq 2, $$ in the last inequality.
Thus, we prove (3.10) for the case $2\|u_0\|_{B^{s}_{p,r}}\geq1$.

Therefore, from the above discussion of the two cases, choosing
$T=\min\{T_1,T_2\}>0$, combining (3.9) and (3.12), we have proved
that $\{u^{(n)}\}_{n\in \mathbb{N}}$ is uniformly bounded in
$C([0,T];B^{s}_{p,r})$. Using Proposition 2.2 (vi) and the fact
$B^{s-2}_{p,r}$ is an algebra as $s>2+\frac{1}{p}$, we have
\begin{eqnarray*}
&&\|[((u^{(n)})^2-(u_x^{(n)})^2)+u^{(n)}]m^{(n+1)}_x+
u_x^{(n)}(m^{(n)})^2+u_x^{(n)}m^{(n)}\|_{B^{s-3}_{p,r}}\\
& \leq& C\{
\|m^{(n+1)}_x\|_{B^{s-3}_{p,r}}[(\|(u^{(n)})^2\|_{B^{s-2}_{p,r}}+\|(u_x^{(n)})^2\|_{B^{s-2}_{p,r}})
+\|u^{(n)}\|_{B^{s-2}_{p,r}}]\\
 &&+
\|u_x^{(n)}\|_{B^{s-3}_{p,r}}(\|(m^{(n)})^2\|_{B^{s-2}_{p,r}}+\|m^{(n)}\|_{B^{s-2}_{p,r}})
\}\\
 & \leq& C[\|u^{(n+1)}\|_{B^{s}_{p,r}}(\|(u^{(n)})\|^2_{B^{s}_{p,r}}+\|u^{(n)}\|_{B^{s}_{p,r}})
 +(\|u^{(n)}\|^3_{B^{s}_{p,r}}+\|u^{(n)}\|^2_{B^{s}_{p,r}})].
 \end{eqnarray*}
From Eq. (3.4), we get $\partial_tm^{(n+1)}\in C\in
([0,T];B^{s-3}_{p,r})$. Hence, $\partial_tu^{(n+1)}\in C\in
([0,T];B^{s-1}_{p,r})$ is uniformly bounded, which yields that the
sequence $\{u^{(n)}\}_{n\in \mathbb{N}}$ is uniformly bounded in
$E^{s}_{p,r}(T)$.

Now it suffices to prove that\\
 $$\{u^{(n)}\}_{n\in \mathbb{N}}\mbox{\ is a
Cauchy sequence in\ } C([0,T];B^{s-1}_{p,r}).$$ Indeed, from Eq.
(3.4), for all $n,l\in \mathbb{N}$, we have
\begin{equation*}
\{\partial_t-[\frac{k_1}{2}((u^{(n+l)})^2-(u_x^{(n+l)})^2)
+\frac{k_2}{2}u^{(n+l)}] \}(m^{(n+l+1)}-m^{(n+1)})=G(t,x),
\end{equation*}
where
$G(t,x):=\{\frac{k_1}{2}[(u^{(n+l)}-u^{(n)})(u^{(n+l)}+u^{(n)})-
(u_x^{(n+l)}-u_x^{(n)})(u_x^{(n+l)}+u_x^{(n)})]
+\frac{k_2}{2}(u^{(n+l)}-u^{(n)})\}\partial_xm^{(n+1)}
+k_1u_x^{(n+l)}(m^{(n+l)}-m^{(n)})(m^{(n+l)}+m^{(n)})+k_1(u_x^{(n+l)}-u_x^{(n)})(m^{(n)})^2
+k_2u_x^{(n+l)}(m^{(n+l)}-m^{(n)})+k_2(u_x^{(n+l)}-u_x^{(n)})m^{(n)}.$

Applying Lemma 2.1 again, for all $t\in [0,T]$, we have
\begin{eqnarray}
&&\|m^{(n+l+1)}-m^{(n+1)}\|_{B^{s-3}_{p,r}}\\ &\leq&
e^{C\int^t_0\|[(u^{(n+l)})^2-(u_x^{(n+l)})^2+u^{(n+l)}](\tau)\|_{B^{s-2}_{p,r}}d\tau}
(\|m_0^{(n+l+1)}-m_0^{(n+1)}\|_{B^{s-3}_{p,r}}\nonumber\\
&&+C\int^t_0e^{-C\int^\tau_0\|[(u^{(n+l)})^2-(u_x^{(n+l)})^2+u^{(n+l)}](\tau')\|_{B^{s-2}_{p,r}}d\tau'}\times
\|G(\tau)\|_{B^{s-3}_{p,r}}d\tau).\nonumber
\end{eqnarray}
Similar to the proof of the estimate of
$\|F(\tau)\|_{B^{s-3}_{p,r}}$ in Lemma 3.1, for
$s>\max\{2+\frac{1}{p},\frac{5}{2},3-\frac{1}{p}\}$ and $s\neq
4+\frac{1}{p}$, we also obtain
\begin{eqnarray*}
\|G(\tau)\|_{B^{s-3}_{p,r}}&\leq&
C\|u^{(n+l)}-u^{(n)}\|_{B^{s-1}_{p,r}}(
\|u^{(n)}\|^2_{B^{s}_{p,r}}+\|u^{(n+1)}\|^2_{B^{s}_{p,r}}\\
&&+\|u^{(n+l)}\|^2_{B^{s}_{p,r}}
+\|u^{(n)}\|_{B^{s}_{p,r}}+\|u^{(n+1)}\|_{B^{s}_{p,r}}+\|u^{(n+l)}\|_{B^{s}_{p,r}}).
\end{eqnarray*}
Inserting the above inequality into (3.13), we have
\begin{eqnarray*}
&&\|u^{(n+l+1)}-u^{(n+1)}\|_{B^{s-1}_{p,r}}\\ &\leq&
e^{C\int^t_0\|[(u^{(n+l)})^2-(u_x^{(n+l)})^2+u^{(n+l)}](\tau)\|_{B^{s-2}_{p,r}}d\tau}
(\|u_0^{(n+l+1)}-u_0^{(n+1)}\|_{B^{s-1}_{p,r}}\nonumber\\
&&+C\int^t_0e^{-C\int^\tau_0\|[(u^{(n+l)})^2-(u_x^{(n+l)})^2+u^{(n+l)}](\tau')\|_{B^{s-2}_{p,r}}d\tau'}
\|u^{(n+l)}-u^{(n)}\|_{B^{s-1}_{p,r}}\\
&&\times( \|u^{(n)}\|^2_{B^{s}_{p,r}}+\|u^{(n+1)}\|^2_{B^{s}_{p,r}}
+\|u^{(n+l)}\|^2_{B^{s}_{p,r}}
+\|u^{(n)}\|_{B^{s}_{p,r}}+\|u^{(n+1)}\|_{B^{s}_{p,r}}\\&&+\|u^{(n+l)}\|_{B^{s}_{p,r}}).\
\end{eqnarray*}
Note that $\{u^{(n)}\}_{n\in \mathbb{N}}$ is uniformly bounded in
$E^{s}_{p,r}(T)$ and
\begin{eqnarray*}
&&\|u_0^{(n+l+1)}-u_0^{(n+1)}\|_{B^{s-1}_{p,r}}\\
&=&\|S_{n+l+1}u_0-S_{n+1}u_0\|_{B^{s-1}_{p,r}}=\|\sum_{q=n+1}^{n+l}\Delta_q
u_0\|_{B^{s-1}_{p,r}}\\
 &\leq&
(\sum_{k\geq-1}2^{k(s-1)r}\|\Delta_k(\sum_{q=n+1}^{n+l}\Delta_q
u_0)\|^r_{L^p})^{\frac{1}{r}}\leq
C(\sum^{n+l+1}_{k=n}2^{-kr}2^{krs}\|\Delta_ku_0\|^r_{L^p})^{\frac{1}{r}}\\
 &\leq& C2^{-n}\|u_0\|_{B^{s}_{p,r}}.
\end{eqnarray*}
Hence, there exists a constant $C_T$ independent of $n,l$ such that
for all $t\in[0,T]$
$$\|(u^{(n+l+1)}-u^{(n+1)})(t)\|_{B^{s-1}_{p,r}}
\leq
C_T(2^{-n}+\int^t_0\|(u^{(n+l)}-u^{(n)})(\tau)\|_{B^{s-1}_{p,r}}d\tau).$$
Arguing by induction with respect to the index $n$, we deduce
\begin{eqnarray*}
&&\|(u^{(n+l+1)}-u^{(n+1)})(t)\|_{B^{s-1}_{p,r}}\\
&\leq&
C_T(2^{-n}\sum_{k=0}^{n}\frac{(2TC_T)^k}{k!}+C_T^{n+1}\int^t_0\frac{(t-\tau)^n}{n!}d\tau)\\
&\leq&(C_T\sum_{k=0}^{n}\frac{(2TC_T)^k}{k!})2^{-n}+C_T\frac{(TC_T)^{n+1}}{(n+1)!},
\end{eqnarray*}
which yields the desired result.

Finally, we can apply the interpolation method, which is similar to
the proof in Lemma 3.1, to the critical case $s=4+\frac{1}{p}$. We
here omit the details. Therefore, we complete the proof of Lemma
3.2.
\end{proof}

Based on the above preparations, we are in position to state the
local existence result of the Cauchy problem (\ref{main}).
\par
\begin{theorem3}
Suppose that $1\leq p,r\leq\infty$,
$s>\max\{2+\frac{1}{p},\frac{5}{2},3-\frac{1}{p}\}$ and $u_0\in
B^{s}_{p,r}$. Then there exists a time $T>0$ such that the Cauchy
problem (\ref{main}) has a unique solution $u\in E^{s}_{p,r}(T)$,
and the mapping $u_0\rightarrow u$ is continuous from $B^{s}_{p,r}$
into
$$C([0,T];B^{s'}_{p,r})\cap C^1([0,T];B^{s'-1}_{p,r})$$
for all $s'<s$ if $r=\infty$ and $s'=s$ otherwise.
\end{theorem3}

\begin{proof}
According to Lemma 3.2, $\{u^{(n)}\}_{n\in \mathbb{N}}\mbox{\ is a
Cauchy sequence in\ } C([0,T];B^{s-1}_{p,r}),$ so it converges to
some function $u\in C([0,T];B^{s-1}_{p,r})$. Thanks to Lemma 3.2 and
Proposition 2.2 (iv) Fatou lemma, we have that $u\in L^\infty
([0,T];B^{s}_{p,r})$. Thus, by the interpolation method, for all
$s'<s,$ we find that $u\in C([0,T];B^{s'}_{p,r})$.

Taking limit in Eq. (3.4), we conclude that $u$ solves Eq.
(\ref{main}) in the sense of $u\in C([0,T];B^{s'-1}_{p,r})$, for all
$s'<s.$ Since $u\in L^\infty ([0,T];B^{s}_{p,r})$ and the fact
$B^{s}_{p,r}$ is an algebra as $s>2+\frac{1}{p}$, the right-hand
side of the following equation
$$m_t-[\frac{k_1}{2}(u^2-u_x^2)+\frac{k_2}{2}u]m_x
= k_1u_xm^2+k_2u_xm,$$ belongs to $L^\infty([0,T];B^{s-2}_{p,r})$.
In particular, for $r<\infty$, Lemma 2.2 enables us to get that
$u\in C([0,T];B^{s'}_{p,r})$ for all $s'\leq s.$ Finally, taking
advantage of Eq. (\ref{main}) again, we obtain that $\partial_t u\in
C([0,T];B^{s-1}_{p,r})$ if $r<\infty$, and in
$L^\infty([0,T];B^{s-1}_{p,r})$ otherwise.

Moreover, the continuity with respect to the initial data in
$$C([0,T];B^{s'}_{p,r})\cap C^1([0,T];B^{s'-1}_{p,r})\quad (\forall s'<s)$$
can be obtained by Lemma 3.1 and a simple interpolation argument.
While the case $s'=s$, a standard of use of a sequence of viscosity
approximate solutions $\{u_\varepsilon\}_{\varepsilon>0}$ for Eq.
(\ref{main}) which converges uniformly in $C([0,T];B^{s}_{p,r})\cap
C^1([0,T];B^{s-1}_{p,r})$ gives the proof of the continuity of
solutions in $E^{s}_{p,r}(T)$. This completes the proof of the
theorem.
\end{proof}

\begin{remark3}
We know that nonhomogeneous Besov spaces contain Sobolev spaces. In
fact, by Fourier-Plancherel formula, we find that the Besov space
$B^s_{2,2}(\R)$ coincides with the Sobolev space $H^s(\R)$.
Therefore, Theorem 3.1 implies that under the assumption $u_0\in
H^s(\R),s>\frac{5}{2}$, we can obtain the local well-posedness
result to Eq. (\ref{main}).
\end{remark3}

\begin{remark3}
The existence time for Eq. (\ref{main}) can be chosen independently
of $s$ in the following sense \cite{Yin}. If
$$u\in C([0,T];H^s)\cap C^1([0,T];H^{s-1})$$
is a solution to Eq. (\ref{main}) with initial data
$H^r,r>\frac{5}{2},r\neq s$, then
$$u\in C([0,T];H^r)\cap C^1([0,T];H^{r-1})$$ with the same time
$T>0$. In particular, if $u_0\in H^\infty$, then $u\in
C([0,T];H^\infty)$.
\end{remark3}

\section{Blow-up scenario}
\newtheorem{theorem4}{Theorem}[section]
\newtheorem{lemma4}{Lemma}[section]
\newtheorem {remark4}{Remark}[section]
\newtheorem{corollary4}{Corollary}[section]
\par
In this section, by using the local well-posedness result of Theorem
3.1 and energy estimates, we present a precise blow-up scenario for
strong solutions to the Cauchy problem (\ref{main}).
\begin{theorem4}
Let $u_0\in H^s(\R), s>\frac{5}{2}$ be given and assume that $T$ is
the maximal existence time of the  solution $u(t,x)$ to Eq.
(\ref{main}) with the initial data $u_0$ guaranteed by Remark 3.1.
When we take $k_1,k_2$ as non-positive constants, then the
corresponding solution $u(t,x)$ blows up in finite time if and only
if
$$\lim\inf_{t\rightarrow T}\inf_{x\in\R}\{mu_x(t,x)\}=-\infty\quad \ \mbox{or}\quad \
\lim\inf_{t\rightarrow T}\inf_{x\in\R}\{u_x(t,x)\}=-\infty.$$
\end{theorem4}

\begin{proof}From Remark 3.1-3.2 and a simple density argument, we
only need to prove that Theorem 4.1 holds true for $s=4$.
Multiplying Eq. (\ref{main}) by $m$, integrating over $\R$ and
integration by parts, we get
\begin{eqnarray}
&&\frac{1}{2}\frac{d}{dt}\int_{\R}m^2dx\\
&=&\frac{k_1}{2}\int_{\R}(u^2-u_x^2)mm_xdx+k_1\int_{\R}u_xm^3dx
+k_2\int_{\R}u_xm^2dx+\frac{k_2}{2}umm_xdx\nonumber\\
&=&-\frac{k_1}{2}\int_{\R}u_xm^3dx+k_1\int_{\R}u_xm^3dx
+k_2\int_{\R}u_xm^2dx-\frac{k_2}{4}\int_{\R}u_xm^2dx\nonumber\\
&=&\frac{k_1}{2}\int_{\R}u_xm^3dx+\frac{3k_2}{4}\int_{\R}u_xm^2dx.\nonumber
\end{eqnarray}
Differentiating Eq. (\ref{main}) with respect to $x$, we deduce
\begin{eqnarray*}
m_{tx}&=&3k_1u_xmm_x-k_1m^3+k_1um^2+\frac{k_1}{2}(u^2-u^2_x)m_{xx}\\
& &+\frac{3k_2}{2}m_xu_x +k_2um-k_2m^2+\frac{k_2}{2}um_{xx}.
\end{eqnarray*}
Multiplying the above equation by $m_x$, and integrating with
respect $x$ over $\R$, we have
\begin{eqnarray}
&&\frac{1}{2}\frac{d}{dt}\int_{\R}m_x^2dx\\
&=&3k_1\int_{\R}u_xmm^2_xdx-k_1\int_{\R}m^3m_xdx+k_1\int_{\R}um^2m_xdx\nonumber\\
&&+\frac{k_1}{2}\int_{\R}(u^2-u^2_x)m_xm_{xx}dx+\frac{3k_2}{2}\int_{\R}m^2_xu_xdx\nonumber\\
&&+ k_2\int_{\R}u
mm_xdx-k_2\int_{\R}m^2m_xdx+\frac{k_2}{2}\int_{\R}um_xm_{xx}dx\nonumber\\
&=&3k_1\int_{\R}u_xmm^2_xdx-\frac{k_1}{3}\int_{\R}u_xm^3dx-\frac{k_1}{4}\int_{\R}m^2_x(u^2-u^2_x)_xdx\nonumber\\
&&+\frac{3k_2}{2}\int_{\R}m^2_xu_xdx-\frac{k_2}{2}\int_{\R}u_xm^2dx-\frac{k_2}{4}\int_{\R}u-xm^2_xdx\nonumber
\end{eqnarray}
\begin{eqnarray*}&=&\frac{5k_1}{2}\int_{\R}u_xmm^2_xdx-\frac{k_1}{3}\int_{\R}u_xm^3dx+\frac{5k_2}{4}\int_{\R}m^2_xu_xdx
-\frac{k_2}{2}\int_{\R}u_xm^2dx.\nonumber
\end{eqnarray*}
From (4.1)-(4.2), we get
\begin{eqnarray*}&&\frac{d}{dt}\int_{\R}(m^2+m_x^2)dx\\
&=&5k_1\int_{\R}u_xmm^2_xdx+\frac{k_1}{3}\int_{\R}u_xm^3dx
+\frac{5k_2}{2}\int_{\R}m^2_xu_xdx+\frac{k_2}{2}\int_{\R}u_xm^2dx.
\end{eqnarray*}
Assume that $T<\infty$ and there exists $N_1,N_2>0$ such that
$mu_x\geq -N_1,u_x\geq -N_2$ for all $(t,x)\in [0,T)\times \R$. Let
us choose $N,k>0$ such that $N:=\max\{N_1,N_2\}$ and
$k:=\max\{-k_1,-k_2\}$. It then follows that
$$\frac{d}{dt}\int_{\R}(m^2+m_x^2)dx\leq 10kN \int_{\R}(m^2+m_x^2)dx.$$
Applying Gronwall's lemma to the above inequality implies for $t\in
[0,T)$,
\begin{eqnarray}
\|m\|^2_{H^1}\leq e^{10NkT} \|m_0\|^2_{H^1}.
\end{eqnarray}
Differentiating Eq. (\ref{main}) with respect to $x$ twice, we
deduce
\begin{eqnarray*}
m_{txx}&=&-6k_1m^2 m_x+5k_1umm_x+4k_1u_xmm_{xx}+3k_1u_xm_x^2\\
&&+k_1u_xm^2+\frac{k_1}{2}(u^2-u^2_x)m_{xxx}
+2k_2u_xm_{xx}-\frac{7k_2}{2}mm_x\\
&&+\frac{5k_2}{2}um_x+k_2u_xm+\frac{k_2}{2}um_{xxx}.
\end{eqnarray*}
Multiplying the above equation by $m_{xx}$, integrating with respect
to $x$ over $\R$, we have
\begin{eqnarray}
&&\frac{1}{2}\frac{d}{dt}\int_{\R}m_{xx}^2dx\\
&=&-6k_1\int_{\R}m^2
m_xm_{xx}dx+5k_1\int_{\R}umm_xm_{xx}dx+4k_1\int_{\R}u_xmm^2_{xx}dx\nonumber\\
&&+3k_1\int_{\R}u_xm_x^2m_{xx}dx
+k_1\int_{\R}u_xm^2m_{xx}dx+\frac{k_1}{2}\int_{\R}(u^2-u^2_x)m_{xx}m_{xxx}dx\nonumber\\
&&+2k_2\int_{\R}u_xm^2_{xx}dx-\frac{7k_2}{2}\int_{\R}mm_xm_{xx}dx
+\frac{5k_2}{2}\int_{\R}um_xm_{xx}dx\nonumber\\
&&+k_2\int_{\R}u_xmm_{xx}dx+\frac{k_2}{2}\int_{\R}um_{xx}m_{xxx}dx\nonumber\\
&=&-6k_1\int_{\R}m^2
m_xm_{xx}dx+5k_1\int_{\R}umm_xm_{xx}dx+4k_1\int_{\R}u_xmm^2_{xx}dx\nonumber\\
&&+k_1\int_{\R}u_xmm^2_{x}dx+2k_1\int_{\R}umm_xm_{xx}dx-k_1\int_{\R}m^2m_xm_{xx}dx\nonumber\\
&&-2k_1\int_{\R}u_xmm_x^2dx+\frac{k_1}{3}\int_{\R}u_xm^3dx-\frac{k_1}{2}\int_{\R}u_xmm^2_{xx}dx\nonumber\\
&&+2k_2\int_{\R}u_xm^2_{xx}dx-\frac{7k_2}{2}\int_{\R}mm_xm_{xx}dx-\frac{5k_2}{4}\int_{\R}u_xm^2_xdx\nonumber
\end{eqnarray}
\begin{eqnarray*}&&-k_2\int_{\R}u_xm^2_xdx+\frac{k_2}{2}\int_{\R}u_xm^2dx-\frac{k_2}{4}\int_{\R}u_xm^2_{xx}dx\nonumber\\
&=&-7k_1\int_{\R}m^2 m_xm_{xx}dx+7k_1\int_{\R}umm_xm_{xx}dx+\frac{7k_1}{2}\int_{\R}u_xmm^2_{xx}dx\nonumber\\
&&-k_1\int_{\R}u_xmm_x^2dx+\frac{k_1}{3}\int_{\R}u_xm^3dx+\frac{7k_2}{4}\int_{\R}u_xm^2_{xx}dx\nonumber\\
&&-\frac{7k_2}{2}\int_{\R}mm_xm_{xx}dx-\frac{9k_2}{4}\int_{\R}u_xm^2dx+\frac{k_2}{2}\int_{\R}u_xm^2dx.\nonumber
\end{eqnarray*}
Combining (4.1)-(4.2) and (4.4), we obtain
\begin{eqnarray*}&&\frac{d}{dt}\int_{\R}(m^2+m_x^2+m^2_{xx})dx\\
&=&-14k_1\int_{\R}m^2 m_xm_{xx}dx+14k_1\int_{\R}umm_xm_{xx}dx+7k_1\int_{\R}u_xmm^2_{xx}dx\\
&&+3k_1\int_{\R}u_xmm_x^2dx+k_1\int_{\R}u_xm^3dx+\frac{7k_2}{2}\int_{\R}u_xm^2_{xx}dx\\
&&-7k_2\int_{\R}mm_xm_{xx}dx-2\int_{\R}u_xm^2dx+\frac{3k_2}{2}\int_{\R}u_xm^2dx.
\end{eqnarray*}
If $mu_x$ and $u_x$ are bounded from below on $[0,T)\times \R$, $i.
e.,$ there exists $N_1,N_2>0$ such that $mu_x\geq -N_1,u_x\geq -N_2$
for all $(t,x)\in [0,T)\times \R$. Similarly, we can choose $N,k>0$
such that $N:=\max\{N_1,N_2\}$ and $k:=\max\{-k_1,-k_2\}$. Then, by
(4.3) and the above equality, we get
\begin{eqnarray*}&&\frac{d}{dt}\int_{\R}(m^2+m_x^2+m^2_{xx})dx\\
&\leq&\frac{21}{2}k
N\int_{\R}(m^2+m_x^2+m^2_{xx})dx+14k(\|m\|^2_{L^\infty}+\|u
m\|_{L^\infty}\\
&&+\|m\|_{L^\infty})\int_{\R}|m_xm_{xx}|dx \\
&\leq&\frac{21}{2}k
N\int_{\R}(m^2+m_x^2+m^2_{xx})dx+14k(\|m\|^2_{H^1}+\|m\|_{H^1})\int_{\R}(m_x^2+m^2_{xx})dx\\
&\leq&
7k[\frac{3N}{2}+2e^{5NkT}\|m_0\|_{H^1}(e^{5NkT}\|m_0\|_{H^1}+1)]\int_{\R}(m^2+m_x^2+m^2_{xx})dx.
\end{eqnarray*}
Hence, applying Gronwall's inequality implies that for all
$t\in[0,T)$
\begin{eqnarray*}
\|u\|^2_{H^4}&\leq& \|m\|^2_{H^2}=\int_{\R}(m^2+m_x^2+m^2_{xx})dx\\
&\leq&
\exp\{7k[\frac{3N}{2}+2e^{5NkT}\|m_0\|_{H^1}(e^{5NkT}\|m_0\|_{H^1}+1)]\}\|m_0\|^2_{H^2}\\
&\leq&
C\exp\{7k[\frac{3N}{2}+2e^{5NkT}\|m_0\|_{H^1}(e^{5NkT}\|m_0\|_{H^1}+1)]\}\|u_0\|^2_{H^4}.
\end{eqnarray*}
The above inequality and Sobolev's embedding theorem ensure that
$u(t,x)$ does not blow up in finite time.

On the other hand, by Sobolev's imbedding theorem, we find that if
$\lim\inf\limits_{t\rightarrow
T}\inf\limits_{x\in\R}\{mu_x(t,x)\}=-\infty\ \mbox{or} \
\lim\inf\limits_{t\rightarrow
T}\inf\limits_{x\in\R}\{u_x(t,x)\}=-\infty,$ then the solution will
blow up in finite time. This completes the proof of the theorem.
\end{proof}

\section{The existence of peaked solutions}
\newtheorem{theorem5}{Theorem}[section]
\newtheorem{lemma5}{Lemma}[section]
\newtheorem {remark5}{Remark}[section]
\newtheorem{corollary5}{Corollary}[section]
\newtheorem{definition5}{Definition}[section]

In order to understand the meaning of a peaked solution to Eq.
(\ref{main}), we first rewrite Eq. (\ref{main}) as
\begin{equation*}
u_t-\frac{k_1}{2}u^2u_x+\frac{k_1}{6}u^3_x-\frac{k_2}{2}uu_x-\frac{k_1}{6}(1-\partial_x^2)^{-1}
u^3_x-\frac{1}{2}\partial_x (1-\partial_x^2)^{-1}
\big(k_1(uu^2_x+\frac{2}{3}u^3)+k_2(u^2+\frac{1}{2}u^2_x)\big)=0.
\end{equation*}
Note that if $p(x)\triangleq \frac{1}{2}e^{-|x|}, x\in\R,$ then
$(1-\partial_x^2)^{-1}f=p\ast f$ for all $f\in L^2.$ From the above
two facts, we can then define the notion of weak solutions as
follows.
\begin{definition5}
Let $u_0\in W^{1,3}$ be given. If $u(t,x)\in L^{\infty}_{loc}([0,T);
W_{loc}^{1,3})$ and satisfies
\begin{eqnarray*}
&&\int_0^T\int_{\R}\Big(u\phi_t-\frac{1}{6}k_1u^3\phi_x-\frac{1}{6}k_1u_x^3\phi
-\frac{1}{4}k_2u^2\phi_x-\frac{1}{2}p\ast
\big(k_1(uu^2_x+\frac{2}{3}u^3)+k_2(u^2+\frac{1}{2}u^2_x)\big)\phi_x\\
&+&\frac{1}{6}k_1 (p\ast u^3_x) \phi \Big)dx
dt+\int_{\R}u_0\phi(0,x)dx=0,
\end{eqnarray*}
for all functions $\phi\in C_c^{\infty}([0,T)\times\R),$  then
$u(t,x)$ is called a weak solution to Eq. (\ref{main}). If $u$ is a
weak solution on $[0,T)$ for every $T>0$, then it is called a global
weak solution.
\end{definition5}

Next, we prove the existence of single peakon to Eq. (\ref{main}).
\begin{theorem5}
The peaked functions of the form
\begin{equation*}
\varphi_c(t,x)=C_1e^{-|x-ct|},
\end{equation*}
where $C_1$ satisfies $\frac{1}{3}k_1C_1^2+\frac{1}{2}k_2C_1+c=0,$
is a global weak solution to Eq. (\ref{main}) in the sense of
Definition 5.1. Moreover, for every time $t\geq 0$, the peaked
solutions $\varphi_c(t,x)$ belongs to $H^1\cap W^{1,\infty}$.
\end{theorem5}
\begin{remark5}
(i)\ For $k_1=0,k_2\neq0,$ we have $C_1=-\frac{2c}{k_2}$. In
particular, if $k_1=0,k_2=-2,$ then we obtain the single peakon
$\varphi_c(t,x)=ce^{-|x-ct|}$ for the CH equation.\\
(ii)\ For $k_1\neq0,$ we easily get
\begin{equation}\label{2.1}
C_1=\frac{-3(\sqrt{3}k_2\pm\sqrt{3k^2_2-16k_1c})}{4\sqrt{3}k_1}.
\end{equation}
If $3k^2_2-16k_1c\geq0$, then the coefficient $C_1$ of the peakons
$\varphi_c$ is a real number. For example, if we choose
$k_1=-2,k_2=0,$ and $c>0$, then we obtain the single peakon
$\varphi_c(t,x)=\pm \sqrt{\frac{3}{2}c}e^{-|x-ct|}$ of the modified
CH equation (\ref{cubic}). If $3k^2_2-16k_1c<0$, then the
coefficient $C_1$ of the peakons $\varphi_c$  is a complex number.
In \cite{Qiao4}, the authors call it as a complex peakon, $i.e.,$
the peakon has the complex coefficient. Thus, we can propose here
the complex peakon for Eq. (\ref{main}), which is not presented in
both the CH equation and the modified CH equation (\ref{cubic}).
\end{remark5}

\begin{proof}
For any test function $\phi(\cdot)\in C_c^{\infty}(\R)$, using
integration by parts, we infer
\begin{eqnarray*}
\int_{\R}e^{-|y|}\phi'(y)dy&=&\int_{-\infty}^0e^y\phi'(y)dy+\int_0^{\infty}e^{-y}\phi'(y)dy\\
&=&e^y\phi(y)\Big|^0_{-\infty}-\int_{-\infty}^0e^y\phi(y)dy+e^{-y}\phi(y)\Big|^{\infty}_{0}
+\int_0^{\infty}e^{-y}\phi(y)dy\\
&=&-\int_{-\infty}^0e^y\phi(y)dy+\int_0^{\infty}e^{-y}\phi(y)dy=\int_{\R}sign(y)e^{-|y|}\phi(y)dy.
\end{eqnarray*}
Thus, for all $t\geq0$, we have
\begin{eqnarray}\label{2.2}
\partial_x\varphi_c(t,x)=-\mbox{sign}(x-ct)\varphi_c(t,x),
\end{eqnarray}
in the sense of distribution $\mathcal{S}'(\R)$. Hence, the peaked
solutions $\varphi_c(t,x)$ belongs to $H^1\cap W^{1,\infty}$. The
same computation as in (\ref{2.2}), for all $t\geq0$, yields,
\begin{eqnarray}\label{2.3}
\partial_t\varphi_c(t,x)=c\ \mbox{sign}(x-ct)\varphi_c(t,x)\in L^\infty.
\end{eqnarray}
If denoting $\varphi_{0,c}(x)\triangleq \varphi_c(0,x)$, then we get
\begin{eqnarray}\label{2.4}
\lim \limits_{t\rightarrow
0^+}\|\varphi_c(t,\cdot)-\varphi_{0,c}(x)\|_{W^{1,\infty}}=0.
\end{eqnarray}
Combining (\ref{2.2})-(\ref{2.4}) and integrating by parts, for
every test function $\phi(t,x)\in C_c^{\infty}([0,\infty)\times
\R)$, we obtain
\begin{eqnarray}\label{2.5}
&&\int_0^{\infty}\int_{\R}\big(\varphi_c\partial_t\phi-\frac{1}{6}k_1\varphi_c^3\partial_x\phi
-\frac{1}{6}k_1(\partial_x\varphi_c)^3\phi-\frac{k_2}{4}\varphi_c^2\partial_x\phi\big)dx
dt+\int_{\R}\varphi_{0,c}(x)\phi(0,x)dx \nonumber\\
&=&-\int_0^{\infty}\int_{\R}\big(\partial_t\varphi_c-\frac{k_1}{2}\varphi^2_c\partial_x\varphi_c
+ \frac{1}{6}k_1(\partial_x\varphi_c)^3
-\frac{k_2}{2}\varphi_c\partial_x\varphi_c\big)\phi dx dt \nonumber\\
&=&-\int_0^{\infty}\int_{\R}\phi\
\mbox{sign}(x-ct)\varphi_c(c+\frac{k_1}{3}\varphi^2_c+\frac{k_2}{2}\varphi_c)dx
dt.
\end{eqnarray}
Form the definition of $\varphi_c$ and $C_1$ satisfying
$\frac{1}{3}k_1C_1^2+\frac{1}{2}k_2C_1+c=0,$ for $ x>ct$, we have
\begin{eqnarray}\label{2.6}
&&\mbox{sign}(x-ct)\varphi_c(c+\frac{k_1}{3}\varphi^2_c+\frac{k_2}{2}\varphi_c)\nonumber\\
&=&C_1e^{-(x-ct)}(c+\frac{k_1}{3}C^2_1e^{-2(x-ct)}+\frac{k_2}{2}C_1e^{-(x-ct)}) \nonumber\\
&=&-\frac{k_1}{3}C_1^3e^{ct-x}-\frac{k_2}{2}C_1^2e^{ct-x}+\frac{k_1}{3}C_1^3e^{3(ct-x)}+\frac{k_2}{2}C_1^2e^{2(ct-x)}.
\end{eqnarray}
Similarly, for $x\leq ct$, we find
\begin{eqnarray}\label{2.7}
&&\mbox{sign}(x-ct)\varphi_c(c+\frac{k_1}{3}\varphi^2_c+\frac{k_2}{2}\varphi_c)\nonumber\\
&=&-C_1e^{x-ct}(c+\frac{k_1}{3}C^2_1e^{2(x-ct)}+\frac{k_2}{2}C_1e^{x-ct}) \nonumber\\
&=&\frac{k_1}{3}C_1^3e^{x-ct}+\frac{k_2}{2}C_1^2e^{x-ct}-\frac{k_1}{3}C_1^3e^{3(x-ct)}-\frac{k_2}{2}C_1^2e^{2(x-ct)}.
\end{eqnarray}

On the other hand, similar to Definition 2.1, we derive
\begin{eqnarray}\label{2.8}
&&-\int_0^{\infty}\int_{\R}\frac{1}{2}(1-\partial^2_x)^{-1}\big(k_1\varphi_c(\partial_x\varphi_c)^2
+\frac{2}{3}k_1\varphi^3_c+k_2\varphi^2_c+\frac{1}{2}k_2(\partial_x\varphi_c)^2\big)\partial_x\phi\nonumber\\
&&+\frac{1}{6}k_1(1-\partial^2_x)^{-1}(\partial_x\varphi_c)^3\phi dx
dt =\int_0^{\infty}\int_{\R}\big[\frac{1}{2}\phi\
\partial_xp\ast\big(k_1\varphi_c(\partial_x\varphi_c)^2
+\frac{k_2}{2}(\partial_x\varphi_c)^2 \big)\nonumber\\
 &&+\phi\ p\ast
\big(k_1\varphi^2_c\partial_x\varphi_c+k_2\varphi_c\partial_x\varphi_c+\frac{k_1}{6}(\partial_x\varphi_c)^3
\big) \big]dx dt.
\end{eqnarray}
From (\ref{2.2}), we have
\begin{eqnarray}\label{2.9}
&&k_1\varphi^2_c\partial_x\varphi_c+k_2\varphi_c\partial_x\varphi_c+\frac{k_1}{6}(\partial_x\varphi_c)^3\nonumber\\
&=&-k_1\mbox{sign}(x-ct)\varphi_c^3-k_2\mbox{sign}(x-ct)\varphi_c^2-\frac{k_1}{6}\mbox{sign}^3(x-ct)\varphi_c^3\nonumber\\
&=&\frac{k_2}{2}\partial_x(\varphi_c^2)+\frac{7}{18}k_1\partial_x(\varphi_c^3).
\end{eqnarray}
Inserting (\ref{2.9}) into (2.8), we obtain
\begin{eqnarray}\label{2.10}
&&-\int_0^{\infty}\int_{\R}\frac{1}{2}(1-\partial^2_x)^{-1}\big(k_1\varphi_c(\partial_x\varphi_c)^2
+\frac{2}{3}k_1\varphi^3_c+k_2\varphi^2_c+\frac{1}{2}k_2(\partial_x\varphi_c)^2\big)\partial_x\phi\nonumber\\
&&+\frac{1}{6}k_1(1-\partial^2_x)^{-1}(\partial_x\varphi_c)^3\phi dx
dt =\int_0^{\infty}\int_{\R}\phi\
\partial_xp\ast\big(\frac{k_1}{2}\varphi_c(\partial_x\varphi_c)^2+\frac{k_2}{4}(\partial_x\varphi_c)^2\nonumber\\
&&+\frac{k_2}{2}\varphi_c^2+\frac{7}{18}k_1\varphi_c^3\big)dx dt.
\end{eqnarray}
Note that $\partial_xp(x)=-\frac{1}{2}\mbox{sign}(x)e^{-|x|},x\in
\R$, we deduce
\begin{eqnarray}\label{2.11}
&&\partial_xp\ast\big(\frac{k_1}{2}\varphi_c(\partial_x\varphi_c)^2+\frac{k_2}{4}(\partial_x\varphi_c)^2
+\frac{k_2}{2}\varphi_c^2+\frac{7}{18}k_1\varphi_c^3\big)(t,x)\nonumber\\
&=&-\frac{1}{2}\int^\infty_{-\infty}\mbox{sign}(x-y)e^{-|x-y|}\big(
(\frac{k_1}{2}\mbox{sign}^2(y-ct)+\frac{7}{18}k_1) C_1^3e^{-3|y-ct|}\nonumber\\
& &+(\frac{k_2}{4}\mbox{sign}^2(y-ct)
+\frac{k_2}{2})C_1^2e^{-2|y-ct|} \big)dy.
\end{eqnarray}
When $x>ct$, we can split the right hand side of (\ref{2.11}) into
the following three parts
\begin{eqnarray*}
&&\partial_xp\ast\big(\frac{k_1}{2}\varphi_c(\partial_x\varphi_c)^2+\frac{k_2}{4}(\partial_x\varphi_c)^2
+\frac{k_2}{2}\varphi_c^2+\frac{7}{18}k_1\varphi_c^3\big)(t,x)\nonumber\\
&=&-\frac{1}{2}\big(\int^{ct}_{-\infty}+\int_{ct}^x+\int_x^\infty\big)\mbox{sign}(x-y)e^{-|x-y|}\big(
(\frac{k_1}{2}\mbox{sign}^2(y-ct)+\frac{7}{18}k_1) C_1^3e^{-3|y-ct|}\nonumber\\
& &+(\frac{k_2}{4}\mbox{sign}^2(y-ct)
+\frac{k_2}{2})C_1^2e^{-2|y-ct|} \big)dy\nonumber\\
 &\triangleq&I_1+I_2+I_3.
\end{eqnarray*}
A direct calculation for each one of the terms $I_i,1\leq i \leq 3,$
yields
\begin{eqnarray*}
I_1&=&-\frac{1}{2}\int^{ct}_{-\infty}e^{-(x-y)}\big(\frac{8}{9}k_1C_1^3e^{3(y-ct)}
+\frac{3}{4}k_2C_1^2e^{2(y-ct)}\big)dy\\
&=&-\frac{4}{9}k_1C_1^3e^{-(x+3ct)}\int^{ct}_{-\infty}e^{4y}dy-
\frac{3}{8}k_2C_1^2e^{-(x+2ct)}\int^{ct}_{-\infty}e^{3y}dy\\
&=&-\frac{k_1}{9}C_1^3e^{ct-x}-\frac{k_2}{8}C_1^2e^{ct-x},
\end{eqnarray*}
\begin{eqnarray*}
I_2&=&-\frac{1}{2}\int^{x}_{ct}e^{-(x-y)}\big(\frac{8}{9}k_1C_1^3e^{-3(y-ct)}
+\frac{3}{4}k_2C_1^2e^{-2(y-ct)}\big)dy\\
&=&-\frac{4}{9}k_1C_1^3e^{-(x-3ct)}\int^{x}_{ct}e^{-2y}dy-
\frac{3}{8}k_2C_1^2e^{-(x-2ct)}\int^{x}_{ct}e^{-y}dy\\
&=&\frac{2}{9}k_1C_1^3(e^{3(ct-x)}-e^{ct-x})+
\frac{3}{8}k_2C_1^2(e^{2(ct-x)}-e^{ct-x}),
\end{eqnarray*}
and
\begin{eqnarray*}
I_3&=&\frac{1}{2}\int^{\infty}_{x}e^{x-y}\big(\frac{8}{9}k_1C_1^3e^{-3(y-ct)}
+\frac{3}{4}k_2C_1^2e^{-2(y-ct)}\big)dy\\
&=&\frac{4}{9}k_1C_1^3e^{x+3ct}\int^{\infty}_{x}e^{-4y}dy+
\frac{3}{8}k_2C_1^2e^{x+2ct}\int^{\infty}_{x}e^{-3y}dy\\
&=&\frac{k_1}{9}C_1^3e^{3(ct-x)}+\frac{k_2}{8}C_1^2e^{2(ct-x)}.
\end{eqnarray*}
By the above equalities $I_1$-$I_3$, for $x>ct$, we have
\begin{eqnarray}\label{2.12}
&&\partial_xp\ast\big(\frac{k_1}{2}\varphi_c(\partial_x\varphi_c)^2+\frac{k_2}{4}(\partial_x\varphi_c)^2
+\frac{k_2}{2}\varphi_c^2+\frac{7}{18}k_1\varphi_c^3\big)(t,x)\nonumber\\
&=&-\frac{k_1}{3}C_1^3e^{ct-x}+\frac{k_1}{3}C_1^3e^{3(ct-x)}-\frac{k_2}{2}C_1^2e^{ct-x}+\frac{k_2}{2}C_1^2e^{2(ct-x)}
\end{eqnarray}
While for the case $x\leq ct$, we can also split the right hand side
of (\ref{2.11}) into the following three parts
\begin{eqnarray*}
&&\partial_xp\ast\big(\frac{k_1}{2}\varphi_c(\partial_x\varphi_c)^2+\frac{k_2}{4}(\partial_x\varphi_c)^2
+\frac{k_2}{2}\varphi_c^2+\frac{7}{18}k_1\varphi_c^3\big)(t,x)\nonumber\\
&=&-\frac{1}{2}\big(\int^{x}_{-\infty}+\int_{x}^{ct}+\int_{ct}^\infty\big)\mbox{sign}(x-y)e^{-|x-y|}\big(
(\frac{k_1}{2}\mbox{sign}^2(y-ct)+\frac{7}{18}k_1) C_1^3e^{-3|y-ct|}\nonumber\\
& &+(\frac{k_2}{4}\mbox{sign}^2(y-ct)
+\frac{k_2}{2})C_1^2e^{-2|y-ct|} \big)dy\nonumber\\
 &\triangleq&II_1+II_2+II_3.
\end{eqnarray*}
We now directly compute each one of the terms $II_i,1\leq i \leq 3,$
as follows
\begin{eqnarray*}
II_1&=&-\frac{1}{2}\int^{x}_{-\infty}e^{-(x-y)}\big(\frac{8}{9}k_1C_1^3e^{3(y-ct)}
+\frac{3}{4}k_2C_1^2e^{2(y-ct)}\big)dy\\
&=&-\frac{4}{9}k_1C_1^3e^{-(x+3ct)}\int^{x}_{-\infty}e^{4y}dy-
\frac{3}{8}k_2C_1^2e^{-(x+2ct)}\int^{x}_{-\infty}e^{3y}dy\\
&=&-\frac{k_1}{9}C_1^3e^{3(x-ct)}-\frac{k_2}{8}C_1^2e^{2(x-ct)},
\end{eqnarray*}
\begin{eqnarray*}
II_2&=&\frac{1}{2}\int^{ct}_{x}e^{x-y}\big(\frac{8}{9}k_1C_1^3e^{3(y-ct)}
+\frac{3}{4}k_2C_1^2e^{2(y-ct)}\big)dy\\
&=&\frac{4}{9}k_1C_1^3e^{x-3ct}\int^{ct}_{x}e^{2y}dy+
\frac{3}{8}k_2C_1^2e^{x-2ct}\int^{ct}_{x}e^{y}dy\\
&=&\frac{2}{9}k_1C_1^3(e^{x-ct}-e^{3(x-ct)})+
\frac{3}{8}k_2C_1^2(e^{x-ct}-e^{2(x-ct)}),
\end{eqnarray*}
and
\begin{eqnarray*}
II_3&=&\frac{1}{2}\int^{\infty}_{ct}e^{x-y}\big(\frac{8}{9}k_1C_1^3e^{-3(y-ct)}
+\frac{3}{4}k_2C_1^2e^{-2(y-ct)}\big)dy\\
&=&\frac{4}{9}k_1C_1^3e^{x+3ct}\int^{\infty}_{ct}e^{-4y}dy+
\frac{3}{8}k_2C_1^2e^{x+2ct}\int^{\infty}_{ct}e^{-3y}dy\\
&=&\frac{k_1}{9}C_1^3e^{x-ct}+\frac{k_2}{8}C_1^2e^{x-ct}.
\end{eqnarray*}
By the above equalities $II_1$-$II_3$, for $x\leq ct$, we obtain
\begin{eqnarray}\label{2.13}
&&\partial_xp\ast\big(\frac{k_1}{2}\varphi_c(\partial_x\varphi_c)^2+\frac{k_2}{4}(\partial_x\varphi_c)^2
+\frac{k_2}{2}\varphi_c^2+\frac{7}{18}k_1\varphi_c^3\big)(t,x)\nonumber\\
&=&\frac{k_1}{3}C_1^3e^{x-ct}-\frac{k_1}{3}C_1^3e^{3(x-ct)}+\frac{k_2}{2}C_1^2e^{x-ct}-\frac{k_2}{2}C_1^2e^{2(x-ct)}.
\end{eqnarray}

Combining (\ref{2.5})-(\ref{2.7}) with (\ref{2.10})-(\ref{2.13}), we
infer that
\begin{eqnarray*}
&&\int_0^{\infty}\int_{\R}\big[\varphi_c\partial_t\phi-\frac{1}{6}k_1\varphi_c^3\partial_x\phi
-\frac{1}{6}k_1(\partial_x\varphi_c)^3\phi-\frac{k_2}{4}\varphi_c^2\partial_x\phi
-\frac{1}{2}(1-\partial^2_x)^{-1}\big(k_1\varphi_c(\partial_x\varphi_c)^2
+\frac{2}{3}k_1\varphi^3_c\\
&&+k_2\varphi^2_c+\frac{1}{2}k_2(\partial_x\varphi_c)^2\big)\partial_x\phi
+\frac{1}{6}k_1(1-\partial^2_x)^{-1}(\partial_x\varphi_c)^3\phi\big]dx
dt+\int_{\R}\varphi_{0,c}(x)\phi(0,x)dx=0
\end{eqnarray*}
for every test function $\phi(t,x)\in C_c^{\infty}([0,\infty)\times
\R)$. This completes the proof of Theorem 5.1.
\end{proof}

\bigskip
\noindent\textbf{Acknowledgements} \ This work was partially
supported by NNSFC (No. 11271382 and No. 10971235), RFDP (No.
200805580014), and the key project of Sun Yat-sen University (No.
c1185).

\end{document}